\documentclass[11pt, a4paper,]{article}
\usepackage{ amsfonts,amsmath}

\newtheorem{prop}{Proposition}[section]

\newcommand{\lig}{{\mathcal{L}}_{\xi}g}

\newcommand{\sdd}{>\!\!\!\triangleleft}

\newcommand{\na}{\stackrel{\star}\nabla}
\newcommand{\nabo}{\stackrel{0}\nabla}
\newcommand{\sn}{\stackrel{\star}{A}_{\xi}}
\newtheorem{defi}{Definition}[section]
\newcommand{\ov}{\overline}
\newtheorem{theo}{Theorem}[section]

\newtheorem{rem}{Remark}[section]

\title{A note on conformal connections on lightlike hypersurfaces.}
\author{C. Atindogbe\footnote{atindogb@iecn.u-nancy.fr,~  Permanent adress: Institut de mathematiques et de sciences Physiques (IMSP), Universit\'e d'Abomey-Calavi(UAC), Benin), 01 BP 613 Porto-Novo, Benin. Email: atincyr@imsp-uac.org} \qquad \qquad L. Berard-Bergery\footnote{berard@iecn.u-nancy.fr}\cr Institut Elie Cartan, Universit\'e Henri Poincar\'e, Nancy I, B.P. 239\cr 54506 Vand\oe uvre-l\`es Nancy Cedex, France}
\date{.}
\begin{document}
\maketitle

\begin{abstract}
\noindent
Degenerate submanifolds of pseudo-Riemannian manifolds are quite difficult to study because there is no prefered connection when the submanifold is not totally geodesic. For the particular case of degenerate totally umbilical hypersurfaces, we show that there are "Weyl" connections adapted to the induced structure on the hypersurface. We begin the study of these with their holonomy.

\end{abstract}


\vspace{0.5cm}

\section{Introduction }
\label{intro}
If $(\ov{M},\ov{g})$ is a Riemannian manifold, any submanifold $M$ of $\ov{M}$ has an induced Riemannian metric $g$. Moreover, the Levi-Civita connection $D$ of the metric $g$ is nothing but the orthogonal projection  along $M$ of the Levi-Civita connection $\ov{D}$ of $\ov{g}$. This is the starting point of the theory of Riemannian submanifolds, and it goes back to Gauss.

If $(\ov{M},\ov{g})$ is a pseudo-Riemannian manifold,the above procedure does not apply as well when the induced "metric" $g$ on $M$ is degenerate. First, there is no longer an orthogonal projection of $T\ov{M}$ onto $TM$ along $M$. Furthermore there is no "Levi-Civita" connection for $g$ in the general case. 

In this paper, we will consider the special case of degenerate hypersurfaces. It is known that $(M,g)$ admits a torsion-free connection $\nabla$ with $g$-parallel ($\nabla g = 0$) if and only if $(M,g)$ is totally geodesic in $(\ov{M},\ov{g})$. When $(M,g)$ is not totally geodesic, one may still choose a torsion-free connection $\nabla$ on $M$ by projecting $\ov{D}$, but there are some choices for such a projection. Such "projected" connections are studied for example in $\cite{DB}$, with the choice of an auxiliary \emph{screen distribution}. 

in this paper we will consider only the special case of totally umbilical hypersurfaces. We remark that a projected connection (as above) cannot be a "Weyl" connection for $(M,g)$. We give a precise definition of Weyl connection in \S~\ref{umbi} (Definition~\ref{conf}), but this is the same as the classical definition for nondegenerate metrics. 

We illustrate that fact in \S~\ref{f3} with the basic example of the light cone at the origin inside the flat $3-$dimensional Lorentz space.

On the other hand, we show in \S~\ref{umb} that there are Weyl connections for any totally umbilical $(M,g)$. We hope that these Weyl connections will be more suited than projected connections for the study of totally umbilical degenerate hypersurfaces. In \S~\ref{hol}, we iluustrate that point of view with a preliminary study of the holonomy of these Weyl connections.

\section{Preliminaries on Lightlike hypersurfaces}
\label{prelim}
Let $M$ be a hypersurface  of an $(n+2)-$dimensional pseudo-Riemannian manifold $(\ov{M},\ov{g})$~of index $0 < \nu < n+2 $. In the classical theory of nondegenerate hypersurfaces, the normal bundle has trivial intersection $\{0\}$ with the tangent one and plays an important role in the introduction of main geometric objects. In case of lightlike (degenerate, null) hypersurfaces, the situation is totally different. The normal bundle $TM^{\perp}$ is a rank-one distribution over $M$: $TM^{\perp}\subset TM$ and then coincide with the so called \emph{radical distribution} $RadTM = TM \cap TM^{\perp}$.  Hence,the induced metric tensor field $g$ is degenerate and has rank $n$. The following characterisation is proved in \cite{DB}.   

\begin{prop}  
\label{prop1}
Let $(M,g)$ be a hypersurface of an $(n+2)-$dimensional pseudo-Riemannian manifold $(\ov{M},\ov{g})$. Then the following assertions are equivalent.
\begin{itemize}
\item[(i)]
$M$ is a lightlike hypersurface of $\ov{M}$.
\item[(ii)]
$g$ has constant rank $n$ on $M$.
\item[(iii)]
$TM^{\perp}= \cup_{x\in M}T_xM^\perp$ is a distribution on $M$.
\end{itemize}
\end{prop}
A complementary  bundle of $TM^\perp$ in $TM$ is a rank $n$ nondegenerate distribution over $M$. It is called a \emph{screen distribution} on $M$  and is often denoted by $S(TM)$. A lightlike hypersurface endowed with a specific screen distribution is denoted by the triple $(M,g,S(TM))$. As $TM^\perp$ lies in the tangent bundle, the following result has an important role in studyng the geometry of a lightlike hypersurface.

\begin{prop}(\cite{DB})
\label{theo1}
Let $(M,g,S(TM))$ be a lightlike hypersurface of $(\ov{M},\ov{g})$ with a given screen distribution $S(TM)$. Then there exists a unique rank $1$ vector subbundle $tr(TM)$ of $T\ov{M}|_{M}$, such that for any non-zero section $\xi$ of $TM^\perp$ on a coordinate neighbourhood ${\mathcal{U}}\subset M$, there exists a unique section $N$ of $tr(TM)$ on ${\mathcal{U}}$ satisfyng 
\begin{equation}
\label{eq1}
\ov{g}(N,\xi)= 1
\end{equation}
and
\begin{equation}
\label{eq2}
\ov{g}(N,N)= \ov{g}(N,W)= 0, \quad \quad \forall~ W \in \Gamma(ST|_{\mathcal{U}}).
\end{equation}
\end{prop}
Here and in the sequel we denote by $\Gamma(E)$ the ${\mathcal{F}}(M)-$module of smooth sections of a vector bundle $E$ over $M$, ${\mathcal{F}}(M)$ being the algebra of smooth functions on $M$. Also, by $\perp$ and $\oplus$ we denote the orthogonal and non-orthogonal direct sum of two vector bundles.  By proposition~\ref{theo1} we may write down the following decompositions.

\begin{equation}
\label{eq3}
TM=S(TM) \perp TM^\perp,
\end{equation}

\begin{equation}
\label{eq4}
T\ov{M}|_{M} = TM \oplus tr(TM)
\end{equation}

and

\begin{equation}
\label{eq4bis}
T\ov{M}|_{M}= S(TM) \perp (TM^\perp \oplus tr(TM))
\end{equation}

As it is well known, we have the following:
\begin{defi}
\label{induced}
Let $(M,g,S(TM))$ be a lightlike hypersurface of $(\ov{M},\ov{g})$ with a given screen distribution $S(TM)$. The induced connection, say $\nabla$, on $M$ is defined by
\begin{equation}
\label{eq48}
\nabla_{X}Y = Q(\ov{\nabla}_XY),
\end{equation} 
where $\ov{\nabla}$ denotes  the Levi-civita connection on $(\ov{M},\ov{g})$ and $Q$ is the projection on $TM$ with respect to the decomposition $(\ref{eq4})$.
\end{defi}

\begin{rem}
\label{rem1}
Notice that the induced connection $\nabla$ on $M$ depends on both $g$ and the specific given screen distribution $S(TM)$ on $M$.
\end{rem}
By respective projections $Q$ and $I-Q$, we have Gauss an Weingarten formulae in the form

\begin{equation}
\label{eq5}
\ov{\nabla}_XY= \nabla_XY + h(X,Y)\qquad \forall X,Y ~ \in \Gamma(TM),
\end{equation}

\begin{equation}
\label{eq6}
\ov{\nabla}_XV= -A_VX + \nabla_X^tV \qquad \forall X ~ \in \Gamma(TM),\quad \forall ~V \in \Gamma(tr(TM)).
\end{equation}
Here, $\nabla_XY$ and $A_{V}X$ belong to $\Gamma(TM)$. Hence

$\bullet$  $h$ is a $\Gamma(tr(TM))$-valued symmetric ${\mathcal{F}}(M)$-bilinear form on $\Gamma(TM)$,

$\bullet$ $A_{V}$ is an  ${\mathcal{F}}(M)$-linear operator on $\Gamma(TM)$, and 

$\bullet$ $\nabla^t$ is a linear connection on the lightlike transversal vector bundle $tr(TM)$. 

Let $P$  denote the projection morphism of $\Gamma(TM)$ on $\Gamma(S(TM))$ with respect to the decomposition (\ref{eq3}). We have 

\begin{equation}
\label{eq9}
\nabla_XPY= \na_XPY + h^{*}(X,PY)\qquad \forall X,Y ~ \in \Gamma(TM),
\end{equation}

\begin{equation}
\label{eq10}
\nabla_X U= -\stackrel{\star}{A}_{U}X + \nabla^{*t}_XU \qquad \forall X ~ \in \Gamma(TM),\quad \forall ~U \in \Gamma(TM^\perp).
\end{equation}

Here $\na_XPY$ and $\stackrel{\star}{A}_{U}X$ belong to $\Gamma(S(TM))$,  $\na$ and $\nabla^{*t}$are linear connection on $S(TM)$ and $TM^\perp$, respectively. Hence

$\bullet$  $h^{*}$ is a $\Gamma(TM^\perp)$-valued  ${\mathcal{F}}(M)$-bilinear form on $\Gamma(TM)\times \Gamma(S(TM))$, and 

$\bullet$ $\stackrel{\star}{A}_{U}$ is a $\Gamma(S(TM))$-valued ${\mathcal{F}}(M)$-linear operator on $\Gamma(TM)$.
 
\noindent
They are the second fundamental form  and the shape operator of the screen distribution, respectively. 

Equivalently, consider a normalizing pair $\{\xi, N\}$ as in the proposition~\ref{theo1}. Then, $(\ref{eq5})$ and $(\ref{eq6})$ take the form

\begin{equation}
\label{eq7}
\ov{\nabla}_XY= \nabla_XY + B(X,Y)N \qquad \forall X,Y ~ \in \Gamma(TM|_{\mathcal{U}}),
\end{equation}
and
\begin{equation}
\label{eq8}
\ov{\nabla}_XN= -A_{N}X+ \tau(X)N \qquad \forall X ~ \in 
\Gamma(TM|_{\mathcal{U}}),
\end{equation}
where we put locally on ${\mathcal{U}}$,
\begin{equation}
\label{eq13}
B(X,Y) = \ov{g}(h(X,Y),\xi)
\end{equation}

\begin{equation}
\label{eq14}
\tau(X) = \ov{g}(\nabla^{t}_XN,\xi)
\end{equation}

\noindent
It is important to stress the fact that the local second fundamental form $B$ in $(\ref{eq13})$ does not depend on the choice of the screen distribution.

We also define (locally) on ${\mathcal{U}}$ the following:
\begin{equation}
\label{eq11}
C(X,PY) = \ov{g}(h^{*}(X,PY),N),
\end{equation}

\begin{equation}
\label{eq12}
\varphi(X) = - \ov{g}(\nabla^{\star t}_X\xi ,N).
\end{equation}

Thus, one has for $~X\in~\Gamma(TM)$
\begin{equation}
\label{eq15}
\nabla_XPY= \na_XPY + C(X,PY)\xi
\end{equation}

\begin{equation}
\label{eq16}
\nabla_X\xi = -\sn X + \varphi(X)\xi 
\end{equation}
It is straighforward to verify that for $~X,Y \in~\Gamma(TM)$
\begin{equation}
\label{eq17}
B(X,\xi) = 0 
\end{equation}

\begin{equation}
\label{eq18}
B(X,Y) = g(\sn X,Y)
\end{equation}

\begin{equation}
\label{eq19}
\sn\xi = 0
\end{equation}

The linear connection $\na$ from (\ref{eq9})is a metric connection on $S(TM)$ and we have for all tangent vector fields $X$, $Y$ and $Z$ in $TM$
\begin{equation}
\label{eq20}
\left(\nabla_{X}g\right)(Y,Z)~=~B(X,Y)\eta(Z) + B(X,Z)\eta(Y).
\end{equation}
with
\begin{equation}
\label{eq49}
\eta(\cdot) = \ov{g}(N,\cdot).
\end{equation}

The induced connection $\nabla$ is torsion-free, but not necessarily $g$-metric. Also, on the geodesibility of $M$ the following is known.

\begin{theo}(\cite[p.88]{DB})
\label{theo2}
Let $(M,g,S(TM))$ be a lightlike hypersurface of a pseudo-Riemannian manifold $(\ov{M},\ov{g} )$. Then the following assertions are equivalent:
\begin{itemize}
\item[(i)] $M$ is totally geodesic.
\item[(ii)]$h$ (or equivalently $B$) vanishes identically on $M$.
\item[(iii)] $\stackrel{\star}{A}_{U}$ vanishes identically on $M$, for any $U~\in \Gamma(TM^\perp)$
\item[(iv)]The  connection $\nabla$ induced by $\ov{\nabla}$ on $M$ is torsion-free and metric.
\item[(v)] $TM^\perp$ is a parallel distribution with respect to $\nabla$.
\item[(vi)]$TM^\perp$ is a Killing distribution on $M$.
\end{itemize}
\end{theo}
It turns out that if $(M,g)$ is not totally geodesic, there is no connection that is, at the same time, torsion-free and $g$-metric. But there is no unicity of such a connection in case there is any.

\section{Some facts about total umbilicity}
\label{umbi}
\begin{defi}
\label{umbil}
A lightlike hypersurface $(M,g)$ is totally umbilical if its second fundamental form $B$ is proportional to the induced metric $g$ pointwise, i.e on each coordinate neighbourhood ${\mathcal{U}}\subset M$, there exists a smooth function $\lambda$ such that 
\begin{equation}
\label{eq21}
B(X,Y)~=~ \lambda g(X,Y), \qquad  \forall ~X,Y\in \Gamma(TM|_{\mathcal{U}}).
\end{equation}
or equivalently,
\begin{equation}
\label{eq22}
\sn PX ~=~\lambda PX, \qquad  \forall ~X \in \Gamma(TM|_{\mathcal{U}}).
\end{equation}

If the function $\lambda$ is nowhere vanishing on $M$, then the latter is said to  be proper totally umbilic.
\end{defi}

It is an easy matter to verify  that this is an intrinsic notion that is independant on ${\mathcal{U}}$, the choice of a screen distribution, $\xi$ (and hence $N$ as in Proposition~\ref{theo1}). 

In some sense, total umbilicity is the nearest situation from being totally geodesic ($\lambda$ is identically zero).
Also, relation (\ref{eq22}) may be written for a given $\xi$ in $Rad(TM)$ as
\begin{equation}
\label{eq23}
\bar{g}(\ov{\nabla}_{X}Y,\xi)=\alpha(\xi)g(X,Y)
\end{equation}
with $\alpha$ a $1$-form on $Rad(TM)$ and $\alpha(\xi)=\lambda$. Hence, $(X,Y)\mapsto \bar{g}(\ov{\nabla}_{X}\xi,Y)= -\alpha(\xi)g(X,Y)$ is a bilinear symmetric form on $\Gamma(TM)$.

\begin{prop}
\label{f1}
A lightlike hypersurface $(M,g)$ is totally umbilic if and only if
there exists a $1$-form $\alpha$ on $RadTM $ such that, for any section $\xi$  of $RadTM$,  
\begin{equation}
\label{eq24}
\lig = -2\alpha(\xi)g.
\end{equation}
Furthermore, $(M,g)$ is proper totally umbilic if and  only if $\alpha$ is everywhere non-zero.
\end{prop}

\noindent
{\bf{Proof.~}} By use of (\ref{eq16}) and (\ref{eq18})we have 

\begin{equation}
\label{eq25}
\lig (X,Y)= -2B(X,Y)
\end{equation}
for all tangent vector fields $X$ and $Y$ in $\Gamma(TM)$.
Then, the equivalence follows (\ref{eq21}).$\square$

In a (Pseudo-) Riemannian setting, manifolds $M^{n}$  with conformal structure $[g]$ and torsion-free connection $D$, such that parallel translation induces conformal transformations, are called Weyl manifolds. In this respect we give the following.

\begin{defi}
\label{conf}
A connection $\nabla$ on a lightlike hypersurface $(M,g)$ is said to be \emph{conformal} if covariant derivative of $g$ is proportional to $g$ in the following precise sense that there exists a $1$-form $\theta$ such that the following,
\begin{equation}
\label{eq26}
\nabla g = \theta \otimes g.
\end{equation}
holds. If in addition, $\nabla$ is torsion-free, it is said to be a \emph{Weyl connection}.
\end{defi}

From Proposition~\ref{f1}, we were tempted to hope that, the induced metric $g$, failing to be parallel with respect to the induced connection, the latter should be conformal. Unfortunately, it is not the case as is shown in the following.
\begin{theo}
\label{f2}
Let $(M,g,S(TM))$ be a lightlike hypersurface with a given screen distribution $S(TM)$. Then, the induced connection is a Weyl connection if and only if  $(M,g)$ is totally geodesic.

In particular, the induced connection on a proper totally umbilical lightlike hypersurface $(M,g,S(TM))$ is never a Weyl connection.
\end{theo}

\noindent
{\bf{Proof.~}}
Assume that the induced connection $\nabla$ satsfies (\ref{eq26}) for some smooth $1$-form $\theta$ defined on $M$. Then, from (\ref{eq20}) we get
\begin{equation}
\label{eq27}
0=\nabla_{X}g(Y,Z)-\theta(X)g(Y,Z)= B(X,Y)\eta(Z) + B(X,Z)\eta(Y) -\theta(X)g(Y,Z)
\end{equation}
for all vector fields $X$,$Y$ and $Z$ in $\Gamma(TM)$.
In (\ref{eq27}), for $Y=\xi \in RadTM$, using $(\ref{eq17})$   leads to $B(X,Z)=0$ for all vector fields $X$ and $Z$ in $\Gamma(TM)$, which is equivalent to saying that $(M,g)$ is totally geodesic. The last assertion is immediate since proper totally umbilic leads to $B$ is everywhere non-zero.$\square$

We shall show that, indeed, there always exists Weyl connections on proper totally umbilic lightlike hypersurfaces. We start by an elementary example.

\section{An elementary example}
\label{f3}
Consider the Lorentz space $\mathbb{R}^{3}_{1}$ that is $\mathbb{R}^{3}$ with the flat Lorentz structure associated with the quadratic form 
\begin{equation}
\label{eq39}
\ov{q}(u) = x^{2} + y^{2}- z^{2},
\end{equation}
for all $u \in \mathbb{R}^{3}$.

Let $\wedge^{2}_{0}$ be the light cone at the origin,
\begin{equation}
\label{eq50}
\wedge^{2}_{0} = \{u \in \mathbb{R}^{3}~ | ~u\ne 0 ~\mbox{and}~ \ov{q}(u) = 0 \}.
\end{equation}
Consider cylindric coordinates system $(r,\theta,z)$ on $\mathbb{R}^{3}$ and the associated frame field $(\partial_{r}, \partial_{\theta},\partial_{z})$ with 
\begin{equation}
\label{eq40}
\left\{
\begin{array}{lcl}
\partial_{r} &=& \cos \theta\, \partial_{x} +  \sin \theta \,\partial_{y}\cr
& & \cr
\partial_{\theta} &=& -\frac{1}{r} \sin \theta \,\partial_{x} +  \frac{1}{r}\cos \theta \,\partial_{y}\cr
& & \cr
\partial_{z}&=& \partial_{z}.
\end{array}
\right.
\end{equation}
Then, it is easy to show that $\xi=\partial_{r} + \partial_{z} $ and $\partial_{\theta} $ are in $T\wedge^{2}_{0}$ and that 
$$ \langle \xi,\xi \rangle = \langle \xi,\partial_{\theta} \rangle = 0$$
where $\ov{g} = \langle~,~\rangle$ denotes the Lorentz metric on $\mathbb{R}_{1}^{3}$. So, $\{ \xi, \partial_{\theta} \}$ is a frame field on $\wedge^{2}_{0}$ and we have $RadT\wedge^{2}_{0}= span\{\xi =\partial_{r} + \partial_{z}   \}$. Thus, $\ov{B}=\{ \xi, \partial_{\theta},\partial_{z} \}$ is a frame field on $\mathbb{R}_{1}^{3} \setminus \{r=0\}$ (which contains $\wedge^{2}_{0}$), with $T\wedge^{2}_{0}= span \{ \xi, \partial_{\theta} \}$. With respect to $\ov{B}$, we have
 
 \begin{equation}
\label{eq44}
\ov{g}_{|{\ov{B}}} = \left(
\begin{array}{ccc}
0& 0& -1 \cr
& &  \cr 
0& \frac{1}{r^{2}}& 0\cr
& &   \cr
-1& 0& -1\cr
\end{array}
\right)
\end{equation}
  
Let $\ov{D}$ denote the Levi-Civita connection on $\mathbb{R}_{1}^{3}$. Then, we have 

\begin{equation}
\label{eq41}
\left\{
\begin{array}{lcl}
\ov{D}_{\xi}\xi &=& 0 \cr
& &\cr 
\ov{D}_{\xi}\partial_{\theta} &=&\displaystyle -\frac{1}{r}\, \partial_{\theta}\cr
& &\cr
\ov{D}_{\partial_{\theta}}\xi &=& \displaystyle -\frac{1}{r}\,\partial_{\theta}\cr
& &\cr
\ov{D}_{\partial_{\theta}}\partial_{\theta}&=& \displaystyle -\frac{1}{r^{3}}\, \partial_{z}.
\end{array}
\right.
\end{equation}
Thus, $\wedge^{2}_{0}$ is proper totally umbilical with $\lambda = \frac{1}{r}$. Now, define on  $\wedge^{2}_{0}$ a connection $\nabo$
as follows.

\begin{equation}
\label{eq42}
\left\{
\begin{array}{lcl}
\nabo{\xi}\xi &=& \displaystyle \frac{2}{r}\, \xi \cr
& &\cr 
\nabo_{\partial_{\theta}}\xi &=& 0\cr
& & \cr
\nabo_{\xi}\partial_{\theta} &=&0\cr
& &\cr
\nabo_{\partial_{\theta}}\partial_{\theta}&=& 0.
\end{array}
\right.
\end{equation}
Let $g$ denote the induced metric on $\wedge^{2}_{0}$.  Then, by straightforward computation using $(\ref{eq42})$ and $(\ref{eq44})$, one gets for all tangent vector fields $X$ and $Y$ in $T\wedge^{2}_{0}$,
\begin{equation}
\label{eq43}
(\nabo_{X}g)(\xi,Y)=0,
\end{equation} 
 and 
 \begin{equation}
\label{eq45}
(\nabo_{X}g)(\theta,\theta)= -2\frac{X\cdot (r)}{r}\cdot \frac{1}{r^{2}}.
\end{equation} 
Combining $(\ref{eq43})$, $(\ref{eq45})$ and $(\ref{eq44})$, we deduce that 
\begin{equation}
\label{eq46}
\nabo g = \omega \otimes g
\end{equation}
with 
\begin{equation}
\label{eq47}
\omega = \frac{2}{r}~\ov{g}~(\cdot~,~ \partial_{z}).
\end{equation}
Thus, since $\nabo$ is torsion-free by construction, it is a Weyl connection on $\wedge^{2}_{0}$.$\square$

From above enumerated facts, one can wonder wether the fact of admitting a Weyl connection characterizes completely proper totally umbilical lightlike hypersurfaces. We are going to give a positive answer to this question.

\section{A total umbilicity criterion}
\label{umb}

\begin{theo}
\label{f4}
For $(M,g)$, lightlike hypersurface of a pseudo-Riemannian manifold $(\ov{M},\ov{g})$ to be to be totally umbilical, it is necessary and sufficient that it admits a Weyl connection. 
\end{theo}

\noindent
{\bf{Proof.~}}
The condition is sufficient. Indeed, suppose there exists a Weyl connection $\nabla$ on $(M,g)$, that is $\nabla$ is torsion-free and there exists on $M$ a smooth $1$-form $\theta$ such that $\nabla~g = \theta \otimes g $. Then, applying this relation to $\xi$, and arbitrary $X$ and $Y$ in $TM$ we get
$$\xi \cdot g(X,Y)- g(\nabla_{\xi}X,Y)-g(X,\nabla_{\xi}Y) = \theta(\xi)g(X,Y)$$
or equivalently, since $\nabla$ is torsion free,
\begin{equation}
\label{eq28}
(\lig)(X,Y)=\theta(\xi)g(X,Y)  + g(\nabla_{X}\xi,Y)+ g(\nabla_{Y}\xi,X).
\end{equation}
Now, we write $\nabla g = \theta \times g$ for $X$, $Y$ in $\Gamma(TM)$ and $\xi \in RadTM$
\begin{eqnarray*}
0  =  \theta(X)g(\xi,Y) &=& (\nabla g)(\xi,Y)\cr
                        &=& X\cdot g (\xi,Y)-g(\nabla_{X}\xi,Y)-g(\xi,\nabla_{X}Y)			 
\end{eqnarray*}
So, $g(\nabla_{X}\xi,Y) = 0$, $\forall~ X,~Y \in TM$, which together with $(\ref{eq28})$ gives $(\lig)(X,Y)= \theta(\xi)g(X,Y)$ which from Proposition~{\ref{f1}} means that $M$ is totally umbilical with $\alpha (\xi)= -\frac{1}{2}\theta(\xi)$.

Conversely, assume $(M,g)$  is a totally umbilical lightlike hypersurface. The $1$-form  $\alpha$ involved in Definition~\ref{umbil} is a section of ${(TM^{\perp})}^{\star}$. Then, the latter is canonically isomorphic to $(T\ov{M}|_{M})/TM$. Also, the projection 
$$T\ov{M}|_{M} \longrightarrow (T\ov{M}|_{M})/TM$$
has contractile fibres. Then, there exists a section $\zeta$ of $T\ov{M}|_{M} $ such that 
$$\forall \xi \in TM^{\perp},~~ \ov{g}(\zeta,\xi ) = \alpha (\xi).$$
Observe that such two sections $\zeta$ differ by exactly one section of $\Gamma(TM)$. Now, let $\theta$ be the smooth $1$-form on $M$ defined by  
\begin{equation}
\label{eq52}
\theta(X)= -2\ov{g}(\zeta, X ),
\end{equation}
and define $\nabla^{\theta}$ as follows.  
\begin{equation}
\label{eq29}
\nabla^{\theta}_{X}Y = \ov{D}_{X}Y -\frac{1}{2}\theta(Y)X-\frac{1}{2}\theta(X)Y- g(X,Y)\zeta
\end{equation}
for all tangent vector fields $X$ and $Y$ in $M$ and $\ov{D}$ the Levi-Civita connection on ambiant space $(\ov{M},\ov{g})$. 

We first show that for $X$ and $Y$ in $TM$, $\nabla^{\theta}_{X}Y \in TM$. Indeed, 
\begin{eqnarray}
\label{eq51}
\ov{g}(\nabla^{\theta}_{X}Y,\xi)  &=& \ov{g}(\ov{D}_{X}Y,\xi)- g(X,Y)\ov{g}(\zeta,\xi)\cr
&\stackrel{(\ref{eq23})}{=}&\alpha(\xi)g(X,Y)
-\alpha(\xi)g(X,Y)=0.
\end{eqnarray}
Thus, $\nabla^{\theta}_{X}Y \in TM$ for $X$ and $Y$ tangent vector fields in $M$. Also, $\nabla^{\theta}$ is clearly a torsion-free connection on $M$. Finally, we show that $\nabla^{\theta}$ is conformal. Let $X$, $Y$ $Z$ in $TM$.We have,
\begin{eqnarray*}
X\cdot g(Y,Z)-g(\nabla^{\theta}_{X}Y,Z)-g(Y,\nabla^{\theta}_{X}Z)&=& \ov{g}(\ov{D}_{X}Y,Z)+\ov{g}(Y,\ov{D}_{X}Z)\cr & & -\ov{g}(\ov{D}_{X}Y,Z)+\frac{1}{2}\theta(Y)\ov{g}(X,Z)\cr & & +\frac{1}{2}\theta(X)\ov{g}(Y,Z)+\ov{g}(X,Y)\ov{g}(\zeta,Z)\cr & &-\ov{g}(Y,\ov{D}_{X}Z)+\frac{1}{2}\theta(Z)\ov{g}(Y,X) \cr & &+ \frac{1}{2}\theta(X)\ov{g}(Y,Z)+\ov{g}(X,Z)\ov{g}(Y,\zeta)
\end{eqnarray*}
that is 
\begin{eqnarray*}
(\nabla^{\theta}_{X}g)(Y,Z)&=&\theta(X)g(Y,Z) +\left[\frac{1}{2}\theta(Z)\ov{g}(Y,X)+\ov{g}(X,Y)\ov{g}(\zeta,Z) \right]\cr & &+\left[\frac{1}{2}\theta(Y)\ov{g}(Z,X)+\ov{g}(X,Z)\ov{g}(\zeta,Y) \right]
\end{eqnarray*}
By $(\ref{eq52})$, the terms in last two brackets are zero.  Hence, $\nabla^{\theta}_{X}g)(Y,Z)= \theta(X)g(Y,Z)$ for all tangent vector fields $X$,$Y$ and  $Z$ in $M$  and the proof is complete.$\square$
\begin{rem}
\label{rem2}
Observe that,
\begin{enumerate}
\item[(a)]if the section $\zeta \in T\ov{M}|_{M}$ is tangent to $M$,we have $\alpha(\xi)$ identically zero and we fall in the totally geodesic case. Otherswise,
\item[(b)] $\zeta$ is nowhere tangent to $M$ and consequently, $\alpha(\xi)$ is everywhere non-zero and $(M,g)$ is proper totally umbilical.
In this case, $S= Ker\theta$ corresponds to the choice of a screen distribution on $M$. For such a screen distribution to be integrable, it is necessary and sufficient that $\theta$ be closed.
\end{enumerate} 
\end{rem}

The above fact gives a good picture of the set of totally umbilical paires $(M,g)$. Let us turn the problem around and relate it to one of the invariants of a given connection on $M$, the holonomy groups \cite{Bes,Joy}. More precisely, let $M$ be a connected  $(n+1)$-manifold, and $\nabla$ a torsion-free  connection on $M$. We ask about necessary and sufficient conditions on $Hol^{0}(\nabla)$ so that $\nabla$ be a Weyl connection of a $1$-degenerate metric $g$ on $M$.

\section{Holonomy groups and $1$-degenerate metrics}
\label{hol}
In the classical case of a manifold with a (non-degenerate) quadratic form $g$ and Weyl connection $\nabla$ with $\nabla g = \theta \otimes g$, the holonomy group of $\nabla$ at point $x$ is included in the conformal group $CO(g_{x}) = \mathbb{R}^{\star}\cdot O(g_{x})$.

Now, let $(M,g)$ be $(n+1)$-manifold $M$ with an everywhere co-rank one metric $g$, and $F$ the frame bundle of $M$. Then, each point of $F$ is $(x,e_{0},\dots,e_{n})$, where $(e_{0},\dots,e_{n})$ is a basis of $T_{x}M$, with $e_{0}$ generating the radical distribution $RadTM$. Now,  define $P$ to be the subset of $F$ for which $(e_{0},\dots,e_{n})$ is a quasi-orthonormal basis with respect to the $1$-degenerate metric $g$ with $sgn(g) = (0,p,n-p)$. Then, $P$ is a $G$-structure on $M$ with fibre 
\begin{equation}
\label{eq30}
G= Aut(TM,g) = \left\{ \left( \begin{matrix}a  &  & ^{t}B \cr
				      0 & &        \cr
				      \vdots & & D \cr            
				      0 & & \end{matrix}  \right), 
a\in \mathbb{R}^{\star},~ B\in \mathbb{R}^{n},~ D \in O(p,n-p)\right\}
\end{equation}
with Lie algebra
\begin{equation}
\label{eq31}
\mathcal{G}= \mathbb{R}^{n}\sdd ~ (\mathbb{R}\oplus so(p,n-p))
\end{equation}

Now, let $(M,\nabla)$ be a connected $(n+1)$-manifold $M$ with a torsion-free connection $\nabla$. We are interested in finding necessary and sufficient conditions on the holonomy group $Hol(\nabla)$ for $\nabla$ to be Weyl connection of a conformal class $[g]$ with $g$ a $1$-degenerate metric on $M$.

Suppose there exists on $M$ a $1$-degenerate metric $g$ and a one form $\theta \in C^{\infty}(T^{\star}M)$  such that $\nabla g = \theta \otimes g$. We first show that the kernel of $g$ that is the radical distribution is $\nabla$-parallel. Let $\xi$ in $RadTM$, $X$, $Y$ in $TM$, we have
\begin{eqnarray*}
g(\nabla_{X}\xi,Y) &=& X\cdot g(\xi,Y)-g(\xi,\nabla_{X}Y)-(\nabla_{X}g)(\xi,Y) \\
&=&-\theta(X)g(\xi,Y)= 0.
\end{eqnarray*}
Then, $RadTM$ is $\nabla$-parallel and it follows similar argument atop of this section that $Hol^{0}(\nabla)$ is a subgroup of $\mathbb{R}^{n}\sdd~\mathbb{R}^{\star}\times (\mathbb{R}^{\star}\cdot SO(n))$, that is any element of $Hol(\nabla)$ has the form
\begin{equation}
\label{eq32}
\left(\begin{tabular}{c|c}
$\alpha$ & $^{t}B $ \cr
\hline
$0$ &   \cr
\vdots & $\lambda~C$\cr
$0$ &                       
		       
		       \end{tabular}\right)
\end{equation}
with $\alpha \in \mathbb{R}^{\star},~B\in \mathbb{R}^{n},~\lambda \in \mathbb{R}^{\star}, C \in O(p,n-p)$.

Conversely, let $\nabla$ be a connection on $M$, and  assume that elements of $Hol(\nabla)$ have the form in (\ref{eq32}). Then, there exists a rank one  subbundle of $TM$, say $L$, that is $\nabla$-parallel. Let 
\[\pi : ~~TM \longrightarrow  TM/L  \]
denote  the natural projection. From $(\ref{eq32})$, the connection $\nabla$ induces on the quotient $TM/L$ a connection $\bar{\nabla}$ explicitly defined by 
\begin{equation}
\label{eq33}
\ov{\nabla}_{X}\ov{Y} = \pi (\nabla_{X}Y)
\end{equation}
for all tangent vectors $X$, $Y$ in $TM$, with $\ov{Y}=\pi(Y)$.
Relation $(\ref{eq33})$ is well defined since for $Y' = Y + U$ with $U \in L$, we have 
\begin{eqnarray*}
\ov{\nabla}_{X}\ov{Y'}& = & \pi(\nabla_{X}Y')  = \pi(\nabla_{X}Y + \nabla_{X}U)  \cr & =& \pi(\nabla_{X}Y) + \pi(\nabla_{X}U) =  \pi(\nabla_{X}Y) = \ov{\nabla}_{X}\ov{Y},
\end{eqnarray*} 
where we used the fact that $L$ is $\nabla$-parallel. Thus, $(\ref{eq33})$ is independant of the choice of  $Y$ in $\ov{Y}$.

On the other hand, it follows $(\ref{eq32})$ and $(\ref{eq33})$ that 
\begin{equation}
\label{eq34}
Hol(\ov{\nabla}) \subset \mathbb{R}^{\star}\cdot O(p,n-p).
\end{equation}

Then, from $(\ref{eq34})$, $\ov{\nabla}$ is a connection on the vector bundle $TM/L$ over $M$, that is compatible with a conformal class $[\ov{g}]$, where $\ov{g}$ is a nondegenerate $(0,2)$  tensor field on the vector bundle  $TM/L$. Alternatively, there exists a one form $\theta \in C^{\infty}(T^{\star}M)$ such that 
\begin{equation}
\label{eq36}
\left(\ov{\nabla}_{X}\ov{g} \right)(\ov{Y},\ov{Z}) = \theta(X)\, \ov{g}(\ov{Y},\ov{Z}) 
\end{equation}
with $\ov{Y}=\pi(Y)$ and $\ov{Z}=\pi(Z)$.

We now define on $M$ a $(0,2)$ tensor $g$ as follows. For tangent vector fields $X$, $Y$ in $TM$, 
\begin{equation}
\label{eq37}
g(X,Y)= \ov{g}(\pi X, \pi Y).
\end{equation}
Clearly, such a $g$ defines on $M$ a degenerate (everywhere) co-rank one metric, with radical distribution  $RadTM = L$. Thus, one has for $X$, $Y$ and $Z$ in $TM$,
\begin{eqnarray}
\label{eq38}
(\nabla_{X}g)(Y,Z)&=&X\cdot g(Y,Z)-g(\nabla_{X}Y,Z)-g(Y,\nabla_{X}Z)\cr
& & \cr
&=&X\cdot \ov{g}(\pi Y,\pi Z)-\ov{g}(\pi(\nabla_{X} Y),\pi Z)-\ov{g}(\pi Y,\pi(\nabla_{X}Z))\cr
& & \cr
&=&X\cdot \ov{g}(\pi Y,\pi Z)-\ov{g}(\ov{\nabla}_{X}\ov{Y},\pi Z)-\ov{g}(\pi Y,\ov{\nabla}_{X}\ov{Z})\cr
& & \cr
&=&(\ov{\nabla}_{X}\ov{g})(\ov{Y},\ov{Z})= \theta(X)\ov{g}(\ov{Y},\ov{Z})\cr
& & \cr
&=&\theta(X)g(Y,Z).
\end{eqnarray}
Thus, $\nabla$ is Weyl connection for the conformal class $[g]$ of the $1$-degenerate metric $g$ on $M$ given by (\ref{eq37}). Thus, we have proved 
\begin{theo}
\label{holo} 
Let $M$ be a connected $(n+1)$-manifold, and $\nabla$ a torsion-free connection on $TM$. Then $\nabla$ is  Weyl connection of a $1$-degenerate metric $g$ on $M$ if and only if $Hol(\nabla)$ is, up to conjugation in $Gl(n+1,\mathbb{R})$, a subgroup of $\mathbb{R}^{n}\sdd ~\left[ \mathbb{R}^{\star} \times (\mathbb{R}^{\star}\cdot O(p,n-p))\right]$.
\end{theo}

\begin{rem}
\label{rem3}
Observe that in theorem~\ref{holo}, one can consider instead the restricted holonomy group $Hol^{0}(\nabla)$. Then, $\nabla$ is  Weyl connection of a $1$-degenerate metric $g$ on $M$ if and only if 

$$ Hol^{0}(\nabla) \subset~ \mathbb{R}^{n}\sdd ~\left[ \mathbb{R}_{+}^{\star} \times (\mathbb{R}_{+}^{\star}\cdot SO_{0}(p,n-p))\right], $$
where $SO_{0}(p,n-p))$ is the connected component of Identity in $SO(p,n-p)$. Recall that, if $M$  is $1$-connected, then $Hol(\nabla) = Hol^{0}(\nabla)$.
\end{rem}

\vspace{1.5cm}
{\bf{Acknowledgments.}}
The first named author (C. Atindogbe) thanks the Agence Universitaire de la Francophonie (AUF) for support with a one year research grant, along with the Institut Elie Cartan (IECN, UHP-Nancy~I) for  research facilities during the completion of this work.

\end{document}